\documentclass [11pt] {article}
\usepackage{amssymb}
\title{Local $\nu$-Euler Derivations and
       \\Deligne's Characteristic Class of
       \\Fedosov Star Products and
       \\Star Products of Special Type}

\author{{\bf Nikolai Neumaier
\thanks{Nikolai.Neumaier@physik.uni-freiburg.de}}\\[3mm]
Fakult\"at f\"ur Physik\\Universit\"at
Freiburg\\Hermann-Herder-Str.3\\79104 Freiburg i.~Br.,
F.~R.~G.\\[3mm]}

\date{FR-THEP-99/3\\[1mm]Revised Version\\[1mm]June 1999\\[6mm]
Dedicated to the memory of\\[1mm]Mosh\'{e} Flato\\[1mm]}

\newcommand {\BEQ} [1] {\begin {equation} \label{#1}}
\newcommand {\EEQ} {\end {equation}}
\newcommand {\BEQAR} [1] {\begin {eqnarray} \label{#1}}
\newcommand {\EEQAR} {\end {eqnarray}}
\newcommand {\Lie} [1] {\mathcal{L}_{#1}}

\newcommand {\id} {{\mathsf {id}}}

\newcommand {\ad} {{\mathrm {ad}}}
\newcommand {\im} {{\mathrm {i}}}
\newcommand {\C} {{\mathrm {C}}}
\newcommand {\Pa} {{\mathrm {P}}}
\newcommand {\starC} {\star_{\mbox{\rm\tiny C}}}
\newcommand {\starP} {\star_{\mbox{\rm\tiny P}}}
\newcommand {\starf} {\star_{\mbox{\rm\tiny flip}}}
\newcommand {\W}   {\mathcal W}
\newcommand {\WL}  {\mbox{$\mathcal W \! \otimes \! \Lambda$}}
\newcommand {\degs} {{\rm deg}_s}
\newcommand {\dega} {{\rm deg}_a}
\newcommand {\degnu} {{\rm deg}_\nu}
\newcommand {\Deg} {{\rm Deg}}
\newenvironment {PROOF}{\small {\sc Proof:}}{{\hspace*{\fill}
                       $\square$}}

\newtheorem {LEMMA} {Lemma} [section]
\newtheorem {PROPOSITION} [LEMMA] {Proposition}
\newtheorem {THEOREM} [LEMMA] {Theorem}
\newtheorem {COROLLARY} [LEMMA] {Corollary}
\newtheorem {DEFINITION}[LEMMA] {Definition}
\newtheorem {REMARK}[LEMMA] {Remark}

\begin{document}
\maketitle
\begin{abstract}
In this paper we explicitly construct local $\nu$-Euler derivations
$\mathsf E_\alpha = \nu \partial_\nu + \Lie{\xi_\alpha} + \mathsf
D_\alpha$, where the $\xi_\alpha$ are local, conformally symplectic
vector fields and the $\mathsf D_\alpha$ are formal series of
locally defined differential operators, for Fedosov star products
on a symplectic manifold $(M,\omega)$ by means of which we are able
to compute Deligne's characteristic class of these star products.
We show that this class is given by
$\frac{1}{\nu}[\omega]+\frac{1}{\nu} [\Omega]$ where $\Omega =
\sum_{i=1}^\infty \nu^i \Omega_i$ is a formal series of closed
two-forms on $M$ the cohomology class of which coincides with the
one introduced by Fedosov to classify his star products. Moreover,
our result implies that the normalisation condition used by Fedosov
does not have any effect on the isomorphy class of the resulting
star product. Finally we consider star products that have
additional algebraic structures and compute the effect of these
structures on the corresponding characteristic classes of these
star products. Specifying the constituents of Fedosov's
construction we obtain star products with these special properties.
Finally we investigate equivalence transformations between such
special star products and prove existence of equivalence
transformations being compatible with the considered algebraic
structures.
\end{abstract}
\tableofcontents
\section{Introduction}
\label{IntroSec}
Since the very beginning of deformation quantisation in the pioneering
articles \cite{BayFla78} by Bayen, Flato, Fronsdal, Lichnerowicz and
Sternheimer there has been not only an immense interest in answering the
question of existence of star products $\star$ (i.~e. formal, associative
deformations of the classical Poisson algebra of complex-valued functions
$C^\infty(M)$ on a symplectic or, more generally, Poisson manifold $M$,
such that in the first order of the formal parameter $\nu$ the commutator
of the star product yields the Poisson bracket) positively, but also in
finding a classification of the star product algebras up to isomorphy of
algebras. Therefore the proofs of existence given by DeWilde and Lecomte
\cite{DeWLec83b}, Fedosov \cite{Fed94,Fed96} in the symplectic case and
recently by Kontsevich \cite{Kon97} in the general case of a Poisson
manifold always contained results on classification. Moreover, there have
been several other results on classification up to equivalence by Nest and
Tsygan \cite{NT95a,NT95b}, Bertelson, Cahen and Gutt \cite{BerCahGut97},
Weinstein and Xu \cite{WX97}.

In the case of deformation quantisations with separation of variables on
K\"{a}hler manifolds Karabegov proved existence and gave a classification using
a deformation of the K\"{a}hler form in \cite{Kar98,Kar99}. Moreover, he has
shown that the Fedosov approach to such star products considered by
Bordemann and Waldmann in \cite{BorWal96a} corresponds to the trivial
deformation in the sense of his classification. This observation gave the
inspiration for our considerations concerning the Deligne class of Fedosov
star products.

In his article \cite{Del95} Deligne has introduced the notions of
\em intrinsic derivation-related \em and \em characteristic class
\em in order to compare the different constructions and
classifications of DeWilde, Lecomte and Fedosov. In his paper
Deligne uses the language of algebraic geometry to approach
deformation theory and proves (cf. \cite[Proposition 3.6.]{Del95})
that the \em relative class \em $c(*)- c(*')$ of two Fedosov star
products being the difference of the characteristic classes of two
Fedosov star products $*$, $*'$ equals $\frac{1}{\nu}(F(*)-
F(*'))$, where $F(*)$ denotes the cohomology class of the
Weyl-curvature Fedosov introduced to classify his star products,
that naturally arises when one constructs a star product using
Fedosov's method.

Recently, Gutt and Rawnsley \cite{GR99} gave an alternative approach to
Deligne's various classes that avoids using methods of algebraic geometry.
Using their methods we succeed in slightly generalising Deligne's result in
proving that Deligne's characteristic class equals $\frac{1}{\nu}$ times
the cohomology class of the Weyl-curvature. We should like to emphasise
that our proof is purely algebraic and does not use any results on sheaf
cohomology except for the de Rham isomorphism relating the second
\v{C}ech cohomology with the second de Rham cohomology. Moreover, we follow
Fedosov's philosophy "Whatever you plan to do on a symplectic manifold $M$,
do it first fibrewise on the tangent space, and pull it then down to $M$ by
means of a compatible Fedosov derivative."

The interest in the relation between the characteristic class and
the Fedosov class is also motivated by the occurrence of the latter
in formulas for canonical traces resp. trace densities obtained by
Halbout in \cite{Hal99} whose results are based on investigations
of invariants in the cyclic cohomology of $M$ made by Connes, Flato
and Sternheimer \cite{CFS92}.

The paper is organised as follows: After a brief summary of
Fedosov's construction of star products on symplectic manifolds we
close Section \ref{FedDelSec} by a short review of the definitions
of Deligne's various classes. Section \ref{EulDerSec} constitutes
the main part of our work, where we give an explicit construction
of local $\nu$-Euler derivations for an arbitrary Fedosov star
product. After these preparations it is an easy task computing
Deligne's derivation-related and characteristic class in Section
\ref{DelClaSec}. As an application of the properties of Deligne's
characteristic class and its relation to Fedosov's Weyl-curvature
we study star products of special type in Section \ref{SpeStaSec}
that satisfy special algebraic identities with respect to complex
conjugation and the mapping $\nu \mapsto -\nu$ changing the sign of
the formal parameter and compute the influence on the corresponding
characteristic classes. Moreover, we can show that there are always
Fedosov star products satisfying these special algebraic identities
the characteristic class of which coincides with a suitably given
element of $\frac{[\omega]}{\nu} + H^2_{\mbox{\rm\tiny
dR}}(M)[[\nu]]$. Considering equivalent star products satisfying
the same algebraic identities with respect to the mappings
mentioned above we can show that there are always equivalence
transformations between these star products commuting with these
mappings. In Appendix \ref{CarSec} we give a short proof of the
deformed Cartan formula that is of great value for our
considerations in Section \ref{DelClaSec}, but seems to be
folklore. A further Appendix \ref{CtwoSec} is added for
completeness giving the computation of the term of the
characteristic class that cannot be determined from the algebraic
considerations in Section \ref{DelClaSec}.

{\bf Conventions:} In what follows we shall use the Einstein summation
convention where the sum over repeated coordinate indices is understood and
the sum ranges from $1$ to $\dim M$. Moreover, we should notice that our
convention for the occurrence of the formal parameter $\nu$ differs from
the one used by Fedosov and in usual physics literature. In our
considerations the star product commutator starts by $\nu$ times the
Poisson bracket, whereas from the viewpoint of physics it is more natural
to have $\im \lambda$ times the Poisson bracket and the formal parameter
$\lambda$ is assumed to be real and directly corresponds to Planck's
constant $\hbar$, so we consider $\nu$ as purely imaginary. The difference
between these conventions is not of great importance as long as one does
not consider star products that have an additional ${}^*$-structure
incorporated by complex conjugation (cf. \cite[Theorem 7.3]{BNW98} for an
example of this situation and Section \ref{SpeStaSec}). In order to
distinguish Fedosov star products from star products in general they shall
be denoted by $*$ whereas the latter ones are denoted by $\star$.
\section{Fedosov Star Products and Deligne's Characteristic
Class}
\label{FedDelSec}
In this rather technical section we shall briefly recall Fedosov's
construction of a star product for a given symplectic manifold
$(M,\omega)$. The notation is mainly the same as in Fedosov's book
\cite{Fed96} and in \cite{Fed94}. In addition we collect the
definitions as they were introduced in \cite{Del95} of \em
Deligne's intrinsic derivation-related class \em and \em Deligne's
characteristic class \em and the relations between them. For proofs
and a detailed discussion of these topics the reader is referred to
the exposition \cite{GR99}.

Let $(M,\omega)$ be a smooth symplectic manifold and define
\BEQ {WeylAlgDef}
    \WL (M) := \left({\mathsf X}_{s=0}^\infty \mathbb C\left(
               \Gamma^\infty \left(\mbox{$\bigvee$}^s T^*M
               \otimes \mbox{$\bigwedge$}
               T^*M\right)\right)\right)[[\nu]] .
\EEQ
If there is no possibility for confusion we simply write $\WL$ and
denote by $\WL^k$ the elements of anti-symmetric degree $k$ and set
$\W := \WL^0$. For two elements $a, b \in \WL$ we define their
pointwise product denoted by $\mu (a \otimes b) = ab$ by the
symmetric $\vee$-product in the first factor and the anti-symmetric
$\wedge$-product in the second factor. Then the degree-maps $\degs$
and $\dega$ with respect to the symmetric and anti-symmetric degree
are derivations of this product. Therefore we shall call $\WL$ a
formally $\mathbb Z \times \mathbb Z$-graded algebra with respect
to the symmetric and anti-symmetric degree. Moreover $(\WL, \mu)$
is super-commutative with respect to the anti-symmetric degree. For
a vector field $X$ we define the symmetric substitution (insertion)
$i_s(X)$ and the anti-symmetric substitution $i_a (X)$ which are
super-derivations of symmetric degree $-1$ resp. $0$ and
anti-symmetric degree $0$ resp. $-1$. Following Fedosov we define
\BEQ {deltaDef}
    \delta := (1 \otimes dx^i) i_s (\partial_i)
    \quad \mbox { and } \quad
    \delta^* := (dx^i \otimes 1) i_a (\partial_i),
\EEQ
where $x^1, \ldots, x^n$ are local coordinates for $M$ and
$\partial_i = \partial_{x^i}$ denotes the corresponding coordinate
vector fields. For $a \in \WL$ with $\degs a = ka$ and $\dega a =
la$ we define
\BEQ {deltaInvDef}
    \delta^{-1}a := \left\{
    \begin {array} {cl}
    \frac{1}{k+l} \delta^*a & \mbox { if } k+l \ne 0 \\
    0 & \mbox { if } k+l=0
    \end {array}
    \right.
\EEQ
and extend $\delta^{-1}$ by linearity. Clearly $\delta^2 = {\delta^*}^2 =
0$. Moreover, we denote by $\sigma: \WL \to C^\infty (M)[[\nu]]$ the
projection onto the part of symmetric and anti-symmetric degree $0$. Then
one has the following `Hodge-decomposition' for any $a \in \WL$ (see e.~g.
\cite [eqn. (2.8)] {Fed94}):
\BEQ {HodgeDecomp}
    a = \delta \delta^{-1} a + \delta^{-1} \delta a + \sigma (a).
\EEQ
Now we consider the fibrewise associative deformation $\circ$ of
the pointwise product having the form
\BEQ {FibProd}
    a \circ b = \mu \circ \exp\left(\frac{\nu}{2}\Lambda^{ij}
    i_s(\partial_i) \otimes i_s(\partial_j)\right)(a \otimes b),
\EEQ
where $\Lambda^{ij}$ denotes the components of the Poisson tensor
corresponding to the symplectic form $\omega$ that are related to
the ones of $\omega$ by the equation $\omega_{kj}\Lambda^{ij}=
\delta^i_k$. Moreover, we define $\dega$-graded super-commutators
with respect to $\circ$ and set $\ad (a) b := [a, b]$.  Now $\degs$
is no longer a derivation of the deformed product $\circ$ but $\Deg
:= \degs + 2 \degnu$ is still a derivation and hence the algebra
$(\WL, \circ)$ is formally $\Deg$-graded where $\degnu := \nu
\partial_\nu$. We shall refer to this degree as total degree.

According to Fedosov's construction of a star product we consider a torsion
free, symplectic connection $\nabla$ on $TM$ that extends in the usual way
to a connection $\nabla$ on $T^*M$ and symmetric resp. anti-symmetric
products thereof. Using this connection we define (using the same symbol as
for the connection) the map $\nabla : \WL \to \WL$ by
\BEQ {Nabladef}
    \nabla:= (1\otimes dx^i) \nabla_{\partial_i},
\EEQ
where $\nabla_{\partial_i}$ denotes the covariant derivative with respect
to $\partial_i$. Then clearly $\nabla$ is globally defined and due to the
property of the connection being symplectic $\nabla$ turns out to be a
super-derivation of anti-symmetric degree $1$ and symmetric and total
degree $0$ of the fibrewise product $\circ$. Moreover $[\delta, \nabla] =
0$ since the connection is torsion free and $\nabla^2 = \frac{1}{2}
[\nabla,\nabla]$ turns out to be an inner super-derivation
\BEQ {Nablaquad}
    \nabla^2 = - \frac{1}{\nu} \ad( R ),
\EEQ
where $R := \frac{1}{4} \omega_{it}R^t_{jkl} dx^i \vee dx^j \otimes dx^k
\wedge dx^l \in \WL^2$ involves the curvature of the connection.
Moreover one has $\delta R = 0 = \nabla R$ as consequences of the
Bianchi identities.

Now remember the following two theorems which are just restatements of
Fedosov's original theorems in \cite[Theorem 3.2, 3.3]{Fed94} resp.
\cite[Theorem 5.3.3]{Fed96}:
\begin {THEOREM} \label{GenFDerivTheo}
Let $\nabla: \WL \to \WL$ be the super-derivation of $\circ$ of
anti-symmetric degree $1$ and total degree $0$ as defined in
(\ref{Nabladef}) satisfying $[\delta,\nabla] = 0$, and $\nabla^2 =
\frac{1}{2} [\nabla, \nabla] = -\frac{1}{\nu} \ad (R)$ with
$R = \frac{1}{4} \omega_{it}R^t_{jkl} dx^i \vee dx^j \otimes dx^k
\wedge dx^l \in \WL^2$ of total degree $2$, satisfying $\delta R =
0 = \nabla R$. Moreover let $\Omega = \sum_{i=1}^\infty
\nu^i\Omega_i$ denote a formal series of closed two-forms on $M$ and let
$s = \sum_{k=3}^\infty s^{(k)} \in \W$ be given with $\sigma(s)=0$ and
$\Deg s^{(k)}= k s^{(k)}$. Then there exists a unique element $r \in
\WL^1$ such that
\BEQ {Genr}
    \delta r = R + \nabla r - \frac{1}{\nu} r \circ r + 1 \otimes
    \Omega \quad \mbox { and } \quad
    \delta^{-1} r = s.
\EEQ
Moreover $r = \sum_{k=2}^\infty r^{(k)}$ with $\Deg r^{(k)} = k r^{(k)}$
satisfies the formula
\BEQ {GenrRecus}
    r = \delta s + \delta^{-1} \left(R + 1 \otimes \Omega +
    \nabla r - \frac{1}{\nu} r \circ r\right)
\EEQ
from which $r$ can be determined recursively. In this case the
Fedosov derivation
\BEQ {GenFedDerivDef}
    \mathcal D := - \delta + \nabla - \frac{1}{\nu} \ad (r)
\EEQ
is a super-derivation of anti-symmetric degree $1$ and has square
zero: $\mathcal D^2 = 0$.
\end {THEOREM}
\begin {THEOREM} \label{GenFTaylorTheo}
Let $\mathcal D = -\delta + \nabla -\frac{1}{\nu}\ad(r): \WL \to
\WL$ be given as in (\ref{GenFedDerivDef}) with $r$ as in
(\ref{Genr}).
\begin {enumerate}
\item Then for any $f \in C^\infty (M)[[\nu]]$ there exists a
      unique element $\tau(f) \in \ker (\mathcal D) \cap \W$ such
      that
      \BEQ {sigmaTaylorL}
          \sigma (\tau(f)) = f
      \EEQ
      and $\tau: C^\infty (M)[[\nu]] \to \ker (\mathcal D)
      \cap\W \subset \mathcal W$
      is $\mathbb C[[\nu]]$-linear and referred to as the
      Fedosov-Taylor series corresponding to $\mathcal D$.

\item In addition we have $\tau(f) = \sum_{k=0}^\infty \tau (f)^{(k)}$ where
      $\Deg\tau (f)^{(k)} = k \tau (f)^{(k)}$ which can
      be obtained by the following recursion formula
      \BEQ {GenTaylorRecurs}
          \begin {array} {c}
              \tau (f)^{(0)} = f \\
              \displaystyle
              \tau (f)^{(k+1)} =
              \delta^{-1} \left(\nabla \tau (f)^{(k)} - \frac{1}{\nu}
              \sum_{l=0}^{k-1}\ad (r^{(l+2)})\tau (f)^{(k-l)}\right).
          \end {array}
      \EEQ

\item Since $\mathcal D$ is a $\circ$-super-derivation of
      anti-symmetric degree $1$ as constructed in Theorem
      \ref {GenFDerivTheo} $\ker \mathcal D \cap \W$ is
      a $\circ$-sub-algebra and a new associative product $*$
      for $C^\infty (M)[[\nu]]$ is defined by pull-back of $\circ$ via
      $\tau$, which turns out to be a star product.
\end {enumerate}
\end {THEOREM}
Observe that in (\ref{Genr}) we allowed for an arbitrary element $s
\in \W$ with $\sigma(s)=0$ that contains no terms of total degree
lower than $3$, as normalisation condition for $r$, i.~e. $\delta^{-1}r=s$
instead of the usually used equation $\delta^{-1}r=0$. In the sequel we
shall especially show that this more general normalisation condition does
not affect the isomorphy class of the resulting star product. In the
following we shall refer to the associative product $*$ defined above as
the Fedosov star product. Moreover, we shall denote by $F(*)$ Fedosov's
characteristic class of the star product $*$ as discussed in \cite[Section
5.3]{Fed96} which is given by $F(*)=[\omega]+ [\Omega]$. There $\omega +
\Omega$ was introduced as the curvature of a certain connection being of
importance for the Fedosov construction, the so-called Weyl-curvature.

Next we collect some basic concepts of characteristic classes for
star products as they can be found in \cite{Del95,GR99}. Deligne's
characteristic class $c(\star)$ of a star product has been
introduced in \cite{Del95} and classifies in a functorial way the
isomorphy classes of star products on a symplectic manifold
$(M,\omega)$. It lies in the affine space $\frac{[\omega]}{\nu} +
H^2_{\mbox{\rm\tiny dR}}(M)[[\nu]]$ and can be calculated by
methods of \v{C}ech cohomology. Let us provide some details of the
calculation as far as they are needed for our purposes. At this
instance we should mention that our conventions, that are as in
\cite{AbrMar85}, differ from those used in \cite{GR99} by a sign in
the Poisson bracket causing the positive sign in front of
$\frac{[\omega]}{\nu}$ in $c(\star)$.

If $\star$ is a star product on the symplectic manifold
$(M,\omega)$ there exists a good open cover $\{
\mathcal U_\alpha\}_{\alpha \in I}$ of $M$ (i.~e. all finite
intersections of the $\mathcal U_\alpha$ are contractible) together
with a family $\{\mathsf E_\alpha\}_{\alpha \in I}$ of local \em
$\nu$-Euler derivations \em of $(C^\infty(\mathcal
U_\alpha)[[\nu]],\star)$ i.~e. a family of derivations $\mathsf
E_\alpha$ of $\star$ over $\mathcal U_\alpha$ having the form
\BEQ {nuEuldef}
    \mathsf E_\alpha = \nu \partial_\nu + \Lie{\xi_\alpha} +
    \mathsf D_\alpha,
\EEQ
where $\xi_\alpha$ is conformally symplectic $(\left.\Lie{
\xi_\alpha} \omega\right|_{\mathcal U_\alpha} = \left.\omega
\right|_{\mathcal U_\alpha})$ and $\mathsf D_\alpha =
\sum_{i = 1}^\infty \nu^i \mathsf D_{\alpha, i}$ is a
formal series of differential operators over $\mathcal U_\alpha$. The
existence of such $\nu$-Euler derivations has already been shown in
\cite{GR99} using cohomological methods, whereas in the case of a Fedosov
star product we are going to give a very direct, purely algebraic proof of
this fact in the next section since for our purposes we need a quite
concrete formula for the differential operators $\mathsf D_\alpha$. As
every $\nu$-linear derivation over a contractible, open set $\mathcal U$ is
of the form $\frac{1}{\nu}\ad_\star ( d )$ with $d \in C^\infty(\mathcal
U)[[\nu]]$ there exist formal functions $d_{\alpha \beta} \in
C^\infty(\mathcal U_\alpha \cap \mathcal U_\beta)[[\nu]]$ fulfilling
\BEQ {Euldiffglei}
    \mathsf E_\alpha - \mathsf E_\beta = \frac{1}{\nu}
    \ad_\star ( d_{\alpha \beta})
\EEQ
over $\mathcal U_\alpha \cap \mathcal U_\beta$. This fact can also be seen
directly from the results of the following two sections. Now, whenever
$\mathcal U_\alpha \cap \mathcal U_\beta \cap
\mathcal U_\gamma \neq \emptyset$ the sums $d_{\alpha \beta
\gamma}= d_{\beta \gamma} - d_{\alpha\gamma}+ d_{\alpha
\beta}$ lie in $\mathbb C [[\nu]]$ and define a $2$-cocycle whose
\v{C}ech class $[d_{\alpha \beta \gamma}] \in
H^2(M, \mathbb C)[[\nu]]$ does not depend on the choices made and
the corresponding class $d(\star)\in H^2_{\mbox{\rm\tiny
dR}}(M)[[\nu]]$ is called \em Deligne's intrinsic
derivation-related class\em.
\begin{DEFINITION}
(cf. \cite[Definition 6.3]{GR99}) Deligne's characteristic class
$c(\star)$ of a star product $\star$ on $(M,\omega)$ is the element
$c(\star) = \frac{[\omega]}{\nu} + \sum_{i=0}^\infty \nu^i
c(\star)^i$ of the affine space $\frac{[\omega]}{\nu} +
H^2_{\mbox{\rm\tiny dR}}(M)[[\nu]]$ defined by
\BEQ {delcladef}
    c(\star)^0 = -2{C_2^-}^\sharp, \qquad \partial_\nu c(\star) =
    \frac{1}{\nu^2} d(\star).
\EEQ
Hereby ${C_2^-}^\sharp$ is the projection onto the second part in
the decomposition $H^2_{\mbox{\rm\tiny Chev, nc}}(C^\infty(M),
C^\infty(M)) = \mathbb C \oplus H^2_{\mbox{\rm\tiny dR}}(M)$ of the
second Chevalley cohomology (null on constants, on $(C^\infty(M),
\{\,{},\,{}\})$ with respect to the adjoint representation) of the
anti-symmetric part $C_2^-(f,g)= \frac{1}{2}\left(C_2(f,g) -
C_2(g,f)\right)$ of the bidifferential operator $C_2$ in the
expansion of $\star$ which is a $2$-cocycle with respect to this
cohomology by the Jacobi-identity for star commutators.
\end{DEFINITION}
\begin{REMARK}
Notice that for Fedosov star products $*$ we have $C_1(f,g)=
\frac{1}{2}\{f,g\}$ implying (cf. \cite[Remark 6.1]{GR99}) that
$C_2^-(f,g) = \rho_2(X_f, X_g)$ for a closed two-form $\rho_2$ on
$M$, where $X_f$ denotes the Hamiltonian vector field with respect
to $\omega$ that corresponds to $f \in C^\infty(M)$, and hence
${C_2^-}^\sharp = [\rho_2]$.
\end{REMARK}

\section{Explicit Construction of Local $\nu$-Euler Derivations}
\label{EulDerSec}
To simplify the notation we use the convention that whenever an
equation contains indices $\alpha, \beta, \gamma$ this means that
it is valid on the intersection of the members of the good open
cover whose indices occur in it. As a first step in the
construction of local $\nu$-Euler derivations we have to find
local, conformally symplectic vector fields $\xi_\alpha$. Since $d
\omega =0$ we can find one-forms $\theta_\alpha$ on each $\mathcal
U_\alpha$ such that $\omega = - d \theta_\alpha$ by the
Poincar\'{e} lemma. Using these local one-forms we can define local
vector fields $\xi_\alpha$ by $i_{\xi_\alpha} \omega = -
\theta_\alpha$. These vector fields obviously satisfy
$\Lie{\xi_\alpha} \omega = \omega$ by Cartan's formula for the Lie
derivative. Using these vector fields we find the following lemma:
\begin{LEMMA}\label{Halphalem}
Let $\mathcal H_\alpha : \WL (\mathcal U_\alpha)
\to \WL (\mathcal U_\alpha)$ be defined by
\BEQ {Halphadef}
    \mathcal H_\alpha := \nu \partial_\nu + \Lie{\xi_\alpha},
\EEQ
then $\mathcal H_\alpha$ is a local (super-)derivation with respect
to the fibrewise product $\circ$ of anti-symmetric and total degree
$0$, i.~e.
\BEQ {Halphaderglei}
    \mathcal H_\alpha (a \circ b) = \mathcal H_\alpha a \circ b +
    a \circ \mathcal H_\alpha b
\EEQ
for all $a, b \in \WL (\mathcal U_\alpha)$. Moreover we have
$[\Lie{\xi_\alpha},\delta]= [\Lie{\xi_\alpha},\delta^*]=0$ and
$[\mathcal H_\alpha, \delta]=[\mathcal H_\alpha, \delta^*]=0$.
\end{LEMMA}
\begin{PROOF}
The proof is a straightforward computation using that $\nu
\partial_\nu$ as well as $\Lie{\xi_\alpha}$ are derivations of the
undeformed product $\mu$ and the equation $\Lie{\xi_\alpha} \Lambda
= - \Lambda$ which follows from $\Lie{\xi_\alpha} \omega = \omega$.
The commutation relations are obvious from the very definitions.
\end{PROOF}

At first sight it might be desirable to construct local derivations
with respect to $*$ by restricting $\mathcal H_\alpha$ to
$C^\infty(\mathcal U_\alpha)[[\nu]]$. In fact this can be done in
some special cases where the connection $\nabla$ is compatible with
the Lie derivative with respect to the vector fields $\xi_\alpha$.
An important example for this situation are homogeneous star
products on cotangent bundles that have been discussed in
\cite{BNW98,BNW99}. But this cannot be done in general since the
failure of the connection to be compatible with the above Lie
derivatives causes that the Fedosov derivation $\mathcal D$ does
not commute with $\mathcal H_\alpha$ and hence $\mathcal H_\alpha$
does not map elements of $\ker(\mathcal D)$ to elements of
$\ker(\mathcal D)$. So we try to extend $\mathcal H_\alpha$ to a
$\circ$-(super-)derivation of anti-symmetric degree $0$ that
commutes with $\mathcal D$. To this end we make the ansatz
\BEQ {Ealphaans}
    \mathcal E_\alpha = \mathcal H_\alpha + \frac{1}{\nu}\ad
    (h_\alpha) = \nu \partial_\nu + \Lie{\xi_\alpha} +
    \frac{1}{\nu}\ad(h_\alpha)
\EEQ
with $h_\alpha \in \W(\mathcal U_\alpha)$ such that $\sigma
(h_\alpha)=0$ and compute $[\mathcal D,\mathcal E_\alpha]$.
\begin{LEMMA}\label{DEalphacomlem}
Let $\mathcal E_\alpha$ be defined as above, then we have
\BEQ {DEalphacomglei}
    [\mathcal D, \mathcal E_\alpha] = \frac{1}{\nu}
    \ad(\mathcal D h_\alpha) + [\nabla, \Lie{\xi_\alpha}] +
    \frac{1}{\nu}\ad(\mathcal H_\alpha r -r).
\EEQ
\end{LEMMA}
\begin{PROOF}
The proof of this formula relies on the fact that $\mathcal D$ is a
super-derivation of anti-symmetric degree $1$ with respect to $\circ$ and
that $\mathcal H_\alpha$ is a (super-)derivation of anti-symmetric degree
$0$ with respect to $\circ$. Moreover we used $[\delta,\mathcal
H_\alpha]=0$ and $[\nabla,\mathcal H_\alpha]= [\nabla, \Lie{\xi_\alpha}]$.
\end{PROOF}

Now we consider the mapping $[\nabla,\Lie{\xi_\alpha}]$ more
closely. The formulas we collect in the following two lemmas are
essential for the whole construction of local $\nu$-Euler
derivations.
\begin{LEMMA}\label{Salphalem1}
For the locally defined vector fields $\xi_\alpha$ the mapping
$[\nabla, \Lie{\xi_\alpha}]$ enjoys the following properties:
\begin{enumerate}
\item
In local coordinates one has
\BEQ {nablaLiecomloc}
    [\nabla,\Lie{\xi_\alpha}] = (dx^j \otimes dx^i)
    i_s((\Lie{\xi_\alpha} \nabla)_{\partial_i}\partial_j)= (dx^j
    \otimes dx^i) i_s(S_\alpha(\partial_i, \partial_j)),
\EEQ
where the local tensor field $S_\alpha \in
\Gamma^\infty(T^*\mathcal U_\alpha \otimes T^*\mathcal U_\alpha
\otimes T\mathcal U_\alpha)$ is defined by
\BEQ {Salphadef}
    S_\alpha(\partial_i,\partial_j)=(\Lie{\xi_\alpha}
    \nabla)_{\partial_i} \partial_j:= \Lie{\xi_\alpha}
    \nabla_{\partial_i}\partial_j - \nabla_{\partial_i}
    \Lie{\xi_\alpha}\partial_j - \nabla_{\Lie{
    \xi_\alpha}\partial_i}\partial_j = R(\xi_\alpha, \partial_i)
    \partial_j + \nabla^{(2)}_{(\partial_i,\partial_j)}\xi_\alpha.
\EEQ
\item
$S_\alpha$ as defined above is symmetric, i.~e. $S_\alpha \in
\Gamma^\infty( \bigvee^2 T^*\mathcal U_\alpha \otimes T\mathcal
U_\alpha)$.
\item
For all $X, Y, Z \in \Gamma^\infty(T \mathcal U_\alpha)$ we have
\BEQ {Salphantisymminsymp}
    \omega(Z, S_\alpha(X,Y))= - \omega(S_\alpha(X,Z),Y).
\EEQ
\end{enumerate}
\end{LEMMA}
\begin{PROOF}
The proof of the local expression for $[\nabla,\Lie{\xi_\alpha}]$ is a
straightforward computation. The last equality in (\ref{Salphadef}) follows
from the torsion freeness of the connection $\nabla$. The fact that
$S_\alpha$ is symmetric is a consequence from the first Bianchi identity
for the connection $\nabla$. Equation (\ref{Salphantisymminsymp}) follows
from a direct computation essentially using $\nabla \omega= 0$ and
$\Lie{\xi_\alpha}\omega =
\omega$.
\end{PROOF}

Now the local tensor fields $S_\alpha$ as defined above naturally
give rise to elements $T_\alpha$ of $\WL(\mathcal U_\alpha)$ of
symmetric degree $2$ and anti-symmetric degree $1$ by
\BEQ {Talphdef}
T_\alpha(Z, Y; X) := \omega(Z, S_\alpha (X,Y)).
\EEQ
In local coordinates this reads $T_\alpha = \frac{1}{2} \omega_{ij}
{S_{\alpha}}_{kl}^j dx^i\vee dx^l \otimes dx^k$, where
${S_{\alpha}}_{kl}^j= dx^j(S_\alpha(\partial_k,\partial_l))$
denotes the components of $S_\alpha$ in local coordinates.
\begin{LEMMA} \label{Salphalem2}
The local tensor field $T_\alpha$ as defined in (\ref{Talphdef})
satisfies the following equations:
\begin{enumerate}
\item
\BEQ {Talphaadgleich}
    \frac{1}{\nu}\ad (T_\alpha) = [\nabla,\Lie{\xi_\alpha}],
\EEQ
\item
\BEQ {Talphagleich}
    T_\alpha = i_a(\xi_\alpha) R + \nabla\left(\frac{1}{2}D
    \theta_\alpha \otimes 1\right),
\EEQ
where the operator of symmetric covariant derivation $D$ is defined by $D:=
dx^i \vee\nabla_{\partial_i}$.
\item
\BEQAR {delnablTalpha}
    \delta T_\alpha = 0 &\textrm{ and }& \nabla T_\alpha =
    \Lie{\xi_\alpha}R - R.
\EEQAR
\end{enumerate}
\end{LEMMA}
\begin{PROOF}
The first assertion easily follows from the properties of $S_\alpha$ given
in Lemma \ref{Salphalem1} by a direct computation. Part ii.) can be easily
proven by direct computation using (\ref{Salphadef}) and the definitions of
$R$ and $T_\alpha$. The equations given in iii.) follow from the
super-Jacobi-identity applied to the equations $[\mathcal
H_\alpha,[\delta,\nabla]]=0$ and $[\mathcal H_\alpha, \frac{1}{2}
[\nabla,\nabla]] = -[\mathcal H_\alpha, \frac{1}{\nu}\ad(R)]$. For the
second equation one has to observe that $R$ does not depend on $\nu$ and
again that $\mathcal H_\alpha$ is a derivation with respect to $\circ$.
Moreover we used the fact that the only central elements of the Fedosov
algebra $\WL$ with respect to $\circ$ with symmetric degree $1$ resp. $2$
are zero.
\end{PROOF}

Collecting our results we have shown that
\BEQ {DEalphcomgleich}
    [\mathcal D, \mathcal E_\alpha] = \frac{1}{\nu}
    \ad(\mathcal D h_\alpha
    +T_\alpha + \mathcal H_\alpha r -r).
\EEQ
Our next aim is to prove that $h_\alpha$ can be chosen such that
$\mathcal D h_\alpha + T_\alpha + \mathcal H_\alpha r - r = 1
\otimes A_\alpha$ where $A_\alpha$ is a formal series of locally
defined one-forms that have to be chosen suitably, since then
$[\mathcal D , \mathcal E_\alpha]=0$. The necessary condition for
this equation to be solvable is $\mathcal D(1 \otimes A_\alpha -
T_\alpha- \mathcal H_\alpha r +r)=0$ since $\mathcal D^2=0$. But
this is also sufficient since the $\mathcal D$-cohomology on
elements $a$ with positive anti-symmetric degree is trivial since
one has the following homotopy formula $\mathcal D
\mathcal D^{-1} a + \mathcal D^{-1} \mathcal Da =a$, where
$\mathcal D^{-1} a := -\delta^{-1}\left(\frac{1}{\id- [\delta^{-1},
\nabla -\frac{1}{\nu}\ad(r)]} a \right)$
(cf. \cite[Theorem 5.2.5]{Fed96}).
\begin{LEMMA}\label{Notbedlem}
Choosing local potentials ${\Theta_i}_\alpha$ for the closed
two-forms $\Omega_i$ on $\mathcal U_\alpha$, and defining
\BEQ {Aalphadef}
    A_\alpha : = (\id - \mathcal H_\alpha)  \Theta_\alpha = (\id -
    \mathcal H_\alpha) \sum_{i=1}^\infty \nu^i {\Theta_i}_\alpha
\EEQ
the equation $\mathcal D(1 \otimes A_\alpha - T_\alpha-
\mathcal H_\alpha r +r)=0$ is fulfilled.
\end{LEMMA}
\begin{PROOF}
Using equation (\ref{Genr}), $[\mathcal H_\alpha, \delta]=0$ and
equation (\ref{Halphaderglei}) as well as Lemma \ref{Salphalem2}
i.),iii.) one computes
\[
\mathcal D \left(\mathcal H_\alpha r - r\right) = 1 \otimes
(\Omega -  \mathcal H_\alpha \Omega) + R - \Lie{\xi_\alpha} R +
[\nabla, \Lie{\xi_\alpha}]r =1 \otimes (\id -  \mathcal H_\alpha
)\Omega - \nabla T_\alpha + \frac{1}{\nu} \ad (T_\alpha) r.
\]
On the other hand we get from $\delta T_\alpha =0$ and $d A_\alpha
= d (\id - \mathcal H_\alpha) \Theta_\alpha = (\id - \mathcal
H_\alpha)\Omega$ that
\[
\mathcal D\left( 1 \otimes A_\alpha - T_\alpha\right) = 1 \otimes
(\id - \mathcal H_\alpha)\Omega - \nabla T_\alpha + \frac{1}{\nu}
\ad(r) T_\alpha,
\]
proving the lemma.
\end{PROOF}

This lemma enables us to prove the following important proposition.
\begin{PROPOSITION} \label{halphaProp}
There are uniquely determined elements $h_\alpha \in
\W(\mathcal U_\alpha)$ such that $\mathcal D h_\alpha = 1
\otimes A_\alpha +r -
\mathcal H_\alpha r - T_\alpha$ and $\sigma (h_\alpha) =0$.
Moreover $h_\alpha$ is explicitly given by
\BEQ {halphagleich}
    h_\alpha = \mathcal D^{-1}\left(  1\otimes A_\alpha + r -
    \mathcal H_\alpha r- T_\alpha \right),
\EEQ
where $\mathcal D^{-1} a = -\delta^{-1}\left(\frac{1}{
\id- [\delta^{-1}, \nabla -\frac{1}{\nu}\ad(r)]}
a \right)$. With these elements $h_\alpha$ the fibrewise, local
$\nu$-Euler derivations $\mathcal E_\alpha = \nu \partial_\nu +
\Lie{\xi_\alpha} + \frac{1}{\nu}\ad(h_\alpha)$ commute with the
Fedosov derivation $\mathcal D$.
\end{PROPOSITION}
\begin{PROOF}
Using the homotopy formula $a= \mathcal D \mathcal D^{-1}a +
\mathcal D^{-1} \mathcal D a$ that is valid for elements $a \in
\WL$ with positive anti-symmetric degree on $1\otimes A_\alpha +
r -\mathcal H_\alpha r- T_\alpha$ we get
\[
1\otimes A_\alpha + r -\mathcal H_\alpha r- T_\alpha = \mathcal D
\mathcal D^{-1}(1\otimes A_\alpha + r -\mathcal
H_\alpha r- T_\alpha)
\]
from the preceding lemma. Since we want the last expression to
equal $\mathcal D h_\alpha$ one gets $h_\alpha =\mathcal
D^{-1}(1\otimes A_\alpha + r -\mathcal H_\alpha - T_\alpha) +
\tau(\varphi_\alpha)$ with arbitrary, locally defined formal
functions $\varphi_\alpha \in C^\infty (\mathcal U_\alpha)[[\nu]]$.
From the demand $\sigma (h_\alpha) = 0$ we get $\sigma (\tau
(\varphi_\alpha))= \varphi_\alpha =0$ since $\mathcal D^{-1}$
raises the symmetric degree and the formula for $h_\alpha$ is
proven. The fact that $\mathcal E_\alpha$ commutes with $\mathcal
D$ now follows from equation (\ref{DEalphcomgleich}).
\end{PROOF}

Using the fibrewise, local $\nu$-Euler derivations $\mathcal
E_\alpha$ we constructed we are in the position to define local
$\nu$-Euler derivations with respect to the Fedosov star product
$*$.
\begin{DEFINITION}
Let $h_\alpha \in \W(\mathcal U_\alpha)$ be given as in equation
(\ref{halphagleich}). Denoting by $\mathcal E_\alpha : \WL
(\mathcal U_\alpha) \to \WL (\mathcal U_\alpha)$ the fibrewise
local $\nu$-Euler derivations $\mathcal E_\alpha = \nu \partial_\nu
+ \Lie{\xi_\alpha} + \frac{1}{\nu}\ad (h_\alpha)$ we define the
mappings $\mathsf E_\alpha : C^\infty(\mathcal U_\alpha)[[\nu]]\to
C^\infty(\mathcal U_\alpha)[[\nu]]$ by
\BEQ {loceudef}
    \mathsf E_\alpha f := \sigma \left( \mathcal E_\alpha \tau (f)
\right)
\EEQ
for $f \in C^\infty (\mathcal U_\alpha)[[\nu]]$.
\end{DEFINITION}
With this definition we get the main result of this section.
\begin{THEOREM}
The mapping $\mathsf E_\alpha$ as defined in equation (\ref{loceudef}) is a
local derivation with respect to the Fedosov star product $*$. Moreover
$\mathsf E_\alpha = \nu \partial_\nu +
\Lie{\xi_\alpha} + \mathsf D_\alpha$, where $\mathsf D_\alpha =
\sum_{i=1}^\infty \nu^i \mathsf D_{\alpha, i}$ is a formal series
of differential operators over $\mathcal U_\alpha$.
\end{THEOREM}
\begin{PROOF}
The fact that $\mathsf E_\alpha$ is a local derivation with respect
to $*$ is obvious from the fact that $\mathcal E_\alpha$ is a local
derivation with respect to $\circ$ and the property of $\mathcal
E_\alpha$ mapping elements in $\W (\mathcal U_\alpha) \cap \ker
\mathcal D$ to elements in $\W (\mathcal U_\alpha) \cap \ker
\mathcal D$ which was achieved by constructing $\mathcal E_\alpha$
such that $[\mathcal D, \mathcal E_\alpha]=0$. The assertion about
the shape of $\mathsf E_\alpha$ follows from the fact that $\sigma$
commutes with $\nu \partial_\nu$ and $\Lie{\xi_\alpha}$ yielding
$\mathsf E_\alpha f = \nu\partial_\nu f + \Lie{\xi_\alpha} f +
\frac{1}{\nu} \sigma \left(\ad (h_\alpha) \tau (f)\right)$. The
fact that the last term involving $h_\alpha$ and $\tau$ defines a
formal series of differential operators is obvious from the
properties of the Fedosov-Taylor series. The only thing one has to
observe is that this formal series starts at order one in the
formal parameter. But this follows from the fact that $h_\alpha$
only contains terms of total degree greater or equal to three,
which is a consequence of $\mathcal D^{-1}$ raising the symmetric
degree, not decreasing the $\nu$-degree and $1\otimes A_\alpha + r
-\mathcal H_\alpha r- T_\alpha$ only containing terms of total
degree greater or equal to two.
\end{PROOF}
\section{Computation of Deligne's Characteristic Class}
\label{DelClaSec}
With the aid of the local $\nu$-Euler derivations we constructed in
the preceding section we are in the position to compute Deligne's
intrinsic derivation-related class $d(*)$ and hence the
characteristic class $c(*)$ for every Fedosov star product $*$ as
defined in Section \ref{FedDelSec}. To this end we have to find
formal functions $d_{\alpha \beta} \in C^\infty(\mathcal
U_\alpha\cap \mathcal U_\beta)[[\nu]]$ such that on $\mathcal
U_\alpha \cap \mathcal U_\beta$ we have $\mathsf E_\alpha
- \mathsf E_\beta = \frac{1}{\nu} \ad_* (d_{\alpha \beta})$.
From the definition of the $\nu$-Euler derivations $\mathsf
E_\alpha$ and the deformed Cartan formula (cf. Appendix
\ref{CarSec}) we have the following:
\begin{LEMMA} \label{EdiffLem}
For $g \in C^\infty(\mathcal U_\alpha\cap \mathcal U_\beta)[[\nu]]$
we have
\BEQ {EdiffGleich}
    \left(\mathsf E_\alpha - \mathsf E_\beta\right) (g) =
    \frac{1}{\nu} \sigma \left(\ad \left( h_\alpha- h_\beta +
    f_{\alpha \beta}+ d f_{\alpha \beta} \otimes 1 +
    \frac{1}{2} D df_{\alpha \beta} \otimes 1 -
    i_a(X_{f_{\alpha \beta}})r\right)\tau(g) \right),
\EEQ
where $f_{\alpha \beta} \in C^\infty( \mathcal U_\alpha\cap
\mathcal U_\beta)$ satisfies $d f_{\alpha \beta} = \theta_\alpha -
\theta_\beta$ and the local one-forms $\theta_\alpha$ satisfy $d
\theta_\alpha = - \omega$.
\end{LEMMA}
\begin{PROOF}
By definition of the $\mathsf E_\alpha$ we have $(\mathsf
E_\alpha-\mathsf E_\beta) (g) = \sigma \left(
\left(\Lie{\xi_\alpha - \xi_\beta} + \frac{1}{\nu}\ad ( h_\alpha -
h_\beta)\right)\tau (g)\right)$. Now on $\mathcal U_\alpha
\cap\mathcal U_\beta$ we have $-d \theta_\alpha = \omega = - d
\theta_\alpha$ and hence by the Poincar\'{e} lemma we can find
locally defined functions $f_{\alpha \beta}$ such that $df_{\alpha
\beta}= \theta_\alpha - \theta_\beta$. Now by definition of the
local vector fields $\xi_\alpha$ we get $d(-f_{\alpha \beta})=
i_{\xi_\alpha- \xi_\beta}\omega$ implying that $\xi_\alpha-
\xi_\beta = X_{-f_{\alpha \beta}}$ is the Hamiltonian vector field
of the function $- f_{\alpha \beta}$. Thus we can apply the
deformed Cartan formula (\ref{LieDerHam}) proven in Proposition
\ref{LieDerProp} and immediately obtain the statement of the lemma
since $\mathcal D \tau(g)=0$ and $i_a(X_{-f_{\alpha
\beta}})\tau(g)=0$ by Theorem \ref{GenFTaylorTheo} and the fact
that $\dega \tau(g)= 0$.
\end{PROOF}

Now we are to show that the term occurring in the argument of $\ad$
in equation (\ref{EdiffGleich}) can be extended by adding a locally
defined formal function $a_{\alpha \beta} \in C^\infty(\mathcal
U_\alpha \cap \mathcal U_\beta)[[\nu]]$ (clearly satisfying
$\ad(a_{\alpha \beta})=0$) such that the whole argument is the
Fedosov-Taylor series $\tau(f_{\alpha \beta}+a_{\alpha
\beta})$ of the local formal function $d_{\alpha \beta}:=f_{\alpha
\beta} + a_{\alpha \beta}$. If we succeed to find such a local
function, equation (\ref{EdiffGleich}) yields $\mathsf E_\alpha -
\mathsf E_\beta = \frac{1}{\nu}\ad_*(d_{\alpha \beta})$ enabling us
to give an expression for Deligne's intrinsic derivation-related
class $d(*)$ of $*$. We thus have to show that $a_{\alpha \beta}$
can be chosen such that
\BEQ {InkerDgleich}
    \mathcal D\left(h_\alpha- h_\beta + a_{\alpha \beta}+
    f_{\alpha \beta}+ d f_{\alpha \beta} \otimes 1 +
    \frac{1}{2} D df_{\alpha \beta} \otimes 1 -
    i_a(X_{f_{\alpha \beta}})r \right)=0.
\EEQ
\begin{LEMMA}\label{Dvonadlem}
With the notations from above we have
\BEQ {Dvonadarg}
    \mathcal D\left(h_\alpha- h_\beta + f_{\alpha \beta}+
    d f_{\alpha \beta} \otimes 1 + \frac{1}{2} D df_{\alpha
    \beta} \otimes 1 - i_a(X_{f_{\alpha \beta}})r \right) =
    1 \otimes ( (i_{\xi_\alpha}\Omega + A_\alpha) -
    ( i_{\xi_\beta}\Omega + A_\beta )),
\EEQ
where $A_\alpha$ is given as in Lemma \ref{Notbedlem}.
\end{LEMMA}
\begin{PROOF}
From the construction of the elements $h_\alpha \in \W (\mathcal
U_\alpha)$ we gave in the preceding section (cf. Proposition
\ref{halphaProp}) we get $\mathcal D \left(h_\alpha -
h_\beta\right)= 1 \otimes (A_\alpha - A_\beta)
- T_\alpha + T_\beta -\Lie{\xi_\alpha -\xi_\beta} r$. Another
straightforward calculation yields $\mathcal D \left( f_{\alpha
\beta} + df_{\alpha \beta}\otimes 1 + \frac{1}{2}D
df_{\alpha \beta}\otimes 1\right)= - \frac{1}{\nu}\ad (r)
\left(df_{\alpha \beta} \otimes 1 + \frac{1}{2}
Ddf_{\alpha \beta}\otimes 1 \right) +
\frac{1}{2}\nabla\left(Ddf_{\alpha\beta}\otimes 1\right)$. Using
the deformed Cartan formula once again combined with equation
(\ref{Genr}) and the definition of $\mathcal D$ we get
\begin{eqnarray*}
\lefteqn{\mathcal D (- i_a(X_{f_{\alpha \beta}})r)=
- (\mathcal D  i_a(X_{f_{\alpha \beta}}) r+
i_a(X_{f_{\alpha \beta}})\mathcal D r) + i_a(X_{f_{\alpha
\beta}})\left(-\delta r+ \nabla r - \frac{1}{\nu}\ad(r)r\right)}\\
&=& \Lie{\xi_\alpha -
\xi_\beta} r + \frac{1}{\nu}\ad\left(-df_{\alpha \beta}\otimes 1
- \frac{1}{2}Ddf_{\alpha \beta}\otimes 1 +
i_a(X_{f_{\alpha \beta}})r \right)r +
i_a(\xi_\alpha-\xi_\beta)\left( R + 1 \otimes \Omega +
\frac{1}{\nu} r \circ r \right)\\
&=&
\Lie{\xi_\alpha -
\xi_\beta} r + \frac{1}{\nu}\ad\left(-df_{\alpha \beta}\otimes 1
- \frac{1}{2}Ddf_{\alpha \beta}\otimes 1\right)r +
i_a(\xi_\alpha-\xi_\beta)R + 1 \otimes i_{\xi_\alpha
-\xi_\beta}\Omega,
\end{eqnarray*}
since $\frac{1}{\nu}i_a(\xi_\alpha - \xi_\beta)(r \circ r) = -
\frac{1}{\nu}\ad(i_a(X_{f_{\alpha \beta}})r)r$. All these results
together with equation (\ref{Talphagleich}) and $df_{\alpha \beta}=
\theta_\alpha - \theta_\beta$ prove the statement of the lemma.
\end{PROOF}

After these preparations we are able to formulate the following
proposition.
\begin{PROPOSITION}
There are locally defined formal functions $d_{\alpha \beta} \in
C^\infty(\mathcal U_\alpha \cap \mathcal U_\beta)[[\nu]]$ such that
$\mathsf E_\alpha - \mathsf E_\beta = \frac{1}{\nu}\ad_*(d_{\alpha
\beta})$. Moreover, these formal functions satisfy $d
d_{\alpha \beta} = d f_{\alpha \beta } + d a_{\alpha \beta} =
\theta_\alpha - \theta_\beta - ((A_\alpha + i_{\xi_\alpha}\Omega )-
(A_\beta + i_{\xi_\beta}\Omega))$. Thus they define a $2$-cocycle
and the image of the corresponding \v{C}ech class under the de Rham
isomorphism, which is just Deligne's intrinsic derivation-related
class, is given by $d(*)= -[\omega]-[ \Omega - \nu \partial_\nu
\Omega]$.
\end{PROPOSITION}
\begin{PROOF}
From Lemma \ref{Dvonadlem} we get that $a_{\alpha \beta}$ has to
satisfy the equation $d a_{\alpha \beta} = - \left( (A_\alpha+
i_{\xi_\alpha} \Omega) - (A_\beta + i_{\xi_\beta}\Omega)\right)$ so
that (\ref{InkerDgleich}) is fulfilled. From the definition of
$A_\alpha$ we get that the right-hand side of this equation is
closed since $d (A_\alpha + i_{\xi_\alpha}\Omega)= \Omega -
\nu\partial_\nu \Omega$. Therefore the existence of $a_{\alpha
\beta}\in C^\infty(\mathcal U_\alpha \cap \mathcal U_\beta)[[\nu]]$
as desired is guaranteed by the Poincar\'{e} lemma. Now we have
$\sigma (h_\alpha- h_\beta + a_{\alpha \beta}+ f_{\alpha
\beta}+ d f_{\alpha \beta} \otimes 1 + \frac{1}{2} D df_{\alpha
\beta} \otimes 1 - i_a(X_{f_{\alpha \beta}})r)=
f_{\alpha \beta} + a_{\alpha \beta}$ and equation
(\ref{InkerDgleich}) is fulfilled implying $(\mathsf
E_\alpha-\mathsf E_\beta)(g) = \frac{1}{\nu}\sigma
(\ad(\tau(f_{\alpha \beta} + a_{\alpha \beta}))\tau(g)=
\frac{1}{\nu}\ad_*(d_{\alpha \beta})g$ by Lemma \ref{EdiffLem}.
The assertion about the corresponding de Rham class is obvious from
the properties of $f_{\alpha \beta}$ and $a_{\alpha \beta}$ we have
already proven namely $d(\theta_\alpha -(A_\alpha +
i_{\xi_\alpha}\Omega)) = -(\omega + \Omega - \nu \partial_\nu
\Omega)$.
\end{PROOF}

From this proposition and from the computation of ${C_2^-}^\sharp$
in Appendix \ref{CtwoSec} we obtain our final result.
\begin{THEOREM} \label{Delclatheo}
Deligne's characteristic class $c(*)$ of a (slightly generalised)
Fedosov star product $*$ as constructed in Section \ref{FedDelSec}
is given by
\BEQ {Delclagleich}
    c(*) = \frac{1}{\nu}[ \omega] + \frac{1}{\nu}[\Omega] =
    \frac{1}{\nu} F(*),
\EEQ
where $F(*)$ denotes Fedosov's characteristic class of the star
product $*$.
\end{THEOREM}
\begin{PROOF}
From the differential equation $\partial_\nu c(*) = \frac{1}{\nu^2}
d(*)$ that relates the derivation-related class to the
characteristic class and from the preceding proposition we get
$c(*)= \frac{1}{\nu}[\omega] + c(*)^0 + \frac{1}{\nu}
\sum_{i=2}^\infty \nu^i [\Omega_i]$. By the result of Proposition
\ref{CtwoProp} we get $c(*)^0 = [\Omega_1]$ proving the theorem.
\end{PROOF}

As an immediate corollary which originally is due to Fedosov (cf.
\cite[Corollary 5.5.4]{Fed96}) we find:
\begin{COROLLARY}
Two Fedosov star products $*$ and $*'$ for $(M,\omega)$ constructed
as in Section \ref{FedDelSec} from the data $(\nabla, \Omega, s)$
and $(\nabla', \Omega', s')$ as in Theorem \ref{GenFDerivTheo} are
equivalent if and only if $[\Omega] = [\Omega']$.
\end{COROLLARY}

\begin{REMARK}
In the case of cotangent bundles $T^*Q$ equipped with the symplectic form
$\omega_0 + \pi^*B_0$, which is a sum of the canonical two-form and the
pull-back of a closed two-form on $Q$, Bordemann, Neumaier, Pflaum and
Waldmann gave a `non-Fedosov' approach to the construction of star products
of a given characteristic class in \cite[Section 4]{BNPW98}. The main idea
of this construction is to start with a star product whose characteristic
class is $[0]$, to apply local equivalence transformations depending on
local potentials of a formal series of closed two-forms $\sum_{i=0}^\infty
\nu^i B_i$ on $Q$ and to show that the local star products obtained fit
together to a globally defined star product.
\end{REMARK}

\section{Star Products of Special Type, their Characteristic
Classes and Equivalence Transformations}
\label{SpeStaSec}In this section we consider star products that
have additional algebraic properties and compute their
characteristic classes showing that these properties give rise to
restrictions on this class. Moreover, we can show that for every
characteristic class satisfying the necessary condition for a star
product of this class to have the desired algebraic properties
there are always Fedosov star products with suitably chosen data
$\Omega, s$ having these properties. Although the following results
might be known they do nevertheless not seem to have appeared in
the literature except for the special case $\Omega=0$ and $s=0$
considered in \cite[Lemma 3.3]{BorWal96a}. In this section $\C: \WL
\to \WL$ shall always denote the complex conjugation, where we
define $\mathrm C \nu := -\nu$ in view of our convention for the
formal parameter being considered as purely imaginary. By $\Pa :
\WL \to \WL$ with $\Pa := (-1)^{\degnu}$ we denote the so-called
$\nu$-parity operator. Using these maps fulfilling $\C^2 = \Pa^2
=\id$ we can define special types of star products:
\begin{DEFINITION}\label{specialdef}
\begin{enumerate}
\item
For a given star product $\star$ for $(M,\omega)$ we define the
star products $\starf$, $\starC$, $\starP$ for $(M, -\omega)$ by
\BEQAR {minusstardef}
    f \starf g &:=& g \star f,\\ f \starC g &:=& \C \left((\C f)
    \star (\C g)\right),\\ f \starP g &:=& \Pa \left((\Pa f)
    \star (\Pa g)\right)= f \star_{-\nu} g =
    \sum_{i=0}^\infty (-\nu)^i C_i(f, g),
\EEQAR
where $f, g \in C^\infty(M)[[\nu]]$ and the bidifferential
operators $C_i$ describe the star product $\star$ by $f \star g =
\sum_{i=0}^\infty \nu^i C_i(f,g)$.
\item
A star product $\star$ is said to have the $\nu$-parity property if
$\Pa$ is an anti-automorphism of $\star$, i.~e.
\BEQ {Nupar}
f \starP g = f \starf g \quad\forall f,g \in C^\infty(M)[[\nu]].
\EEQ
\item
A star product $\star$ is said to have a ${}^*$-structure
incorporated by complex conjugation if $\C$ is an anti-automorphism
of $\star$, i.~e.
\BEQ {Sternstr}
f \starC g = f \starf g \quad\forall f,g \in C^\infty(M)[[\nu]].
\EEQ
\item
A star product $\star$ is called of Weyl type if it has the
$\nu$-parity property and has a ${}^*$-structure incorporated by
complex conjugation.
\end{enumerate}
\end{DEFINITION}

Using these definitions we find:
\begin{LEMMA} \label{classlem}
\begin{enumerate}
\item
The characteristic classes of $\starf$, $\starC$, $\starP$ are
related to the characteristic class $c(\star)$ of $\star$ by the
following equations:
\BEQAR {minusstarclass}
    c(\starf) &=& - c(\star),\\ c(\starC) &=& \C c(\star),\\
    c(\starP)(\nu)&=&c(\star_{-\nu})(\nu) = c(\star)(-\nu) = \Pa
    \left(c(\star)(\nu)\right).
\EEQAR
\item
The characteristic class of a star product $\star$ that has the
$\nu$-parity property satisfies
\BEQ {Nuparclass}
    \Pa c(\star) = - c(\star),
\EEQ
and hence $c(\star) = \frac{[\omega]}{\nu} +
\sum_{l=0}^\infty \nu^{2l+1}c(\star)^{2l+1}$, i.~e.
$c(\star)^{2l}=[0]$ for all $l \in \mathbb N$.
\item
The characteristic class of a star product $\star$ that has $\C$ as
${}^*$-structure satisfies
\BEQ {Sternstrclass}
    \C c(\star) = - c(\star),
\EEQ
and hence $c(\star)^{2l} = - \C c(\star)^{2l}$ and
$c(\star)^{2l+1}= \C c (\star)^{2l+1}$ for all $l \in \mathbb N$.
\item
The characteristic class of a star product $\star$ that is of Weyl
type satisfies
\BEQ {Weylclass}
    \Pa c(\star) = - c(\star) \quad \textrm{ and }\quad
    \C c(\star) = - c(\star),
\EEQ
and hence $c(\star)^{2l}=[0]$ and $c(\star)^{2l+1}=
\C c(\star)^{2l+1}$ for all $l \in \mathbb N$.
\end{enumerate}
\end{LEMMA}
\begin{PROOF}
The proof of part i.) relies on the observation, that local
$\nu$-Euler derivations $\mathsf E_\alpha$ of $\star$ yield such
derivations for $\starf$, $\starC$ and $\starP$ given by $\mathsf
E_\alpha$, $\C \mathsf E_\alpha \C$ and $\Pa \mathsf E_\alpha
\Pa$. With these derivations one easily finds $d(\starf)= -
d(\star)$, $d(\starC) = - \C d(\star)$ and $d(\starP) = -
\Pa d(\star)$. From the definition of the characteristic class
relating the derivation-related class $d$ with $c$ and the obvious
observations that $c(\starf)^0 = - c(\star)^0$, $c(\starC)^0 = \C
c(\star)^0$ and $c(\starP)^0 = c(\star)^0$ one gets the asserted
statements. The assertions ii.),iii.) and iv.) are obvious from
part i.) and Definition \ref{specialdef} ii.),iii.) and iv.).
\end{PROOF}

The preceding lemma states that in general there are equivalence
classes of star products corresponding to the characteristic
classes $c(\star)$ that contain no representatives (i.~e. star
products $\star$ with this characteristic class) satisfying the
conditions (\ref{Nupar}) resp. (\ref{Sternstr}), namely those whose
characteristic classes do not satisfy the equations
(\ref{Nuparclass}) resp. (\ref{Sternstrclass}). Vice versa the
following proposition states that for every class $c \in
\frac{[\omega]}{\nu}+H^2_{\mbox{\rm\tiny dR}}(M)[[\nu]]$ enjoying
the properties $\C c = - c$ resp. $\Pa c = -c$ one can find even
Fedosov star products having the characteristic class $c$ and
satisfying the conditions (\ref{Sternstr}) resp. (\ref{Nupar}).
\begin{PROPOSITION} \label{Fedclassprop}
\begin{enumerate}
\item
For all $c \in \frac{[\omega]}{\nu}+H^2_{\mbox{\rm\tiny
dR}}(M)[[\nu]]$ with $\Pa c = -c$ there are Fedosov star products
$*$ for $(M, \omega)$ with
\[
c(*) = c \quad \textrm{ and } \quad \Pa \left( (\Pa f) * (\Pa
g)\right) = g * f \quad \textrm{ for all }f, g \in
C^\infty(M)[[\nu]].
\]
\item
For all $c \in \frac{[\omega]}{\nu}+H^2_{\mbox{\rm\tiny
dR}}(M)[[\nu]]$ with $\C c = - c$ there are Fedosov star products
$*$ for $(M,\omega)$ with
\[
c(*) = c \quad \textrm{ and } \quad \C \left( (\C f )* (\C g)
\right) = g * f \quad \textrm{ for all }f, g
\in C^\infty(M)[[\nu]].
\]
\item
For all $c \in \frac{[\omega]}{\nu}+H^2_{\mbox{\rm\tiny
dR}}(M)[[\nu]]$ with $\Pa c = - c = \C c$ there are Fedosov star
products $*$ for $(M,\omega)$ with
\[
c(*)=c\quad \textrm{ and }\quad \Pa \left( (\Pa f) * (\Pa g)\right)
= g * f =\C \left( (\C f )* (\C g)\right)\quad \textrm{ for all }
f, g \in C^\infty(M)[[\nu]].
\]
\end{enumerate}
\end{PROPOSITION}
\begin{PROOF}
For the proof we first observe that the fibrewise product $\circ$
satisfies $\C ( (\C a)\circ (\C b)) = \Pa ( (\Pa a) \circ (\Pa b))=
(-1)^{k l} b \circ a$ for all $a, b \in \WL$ with $\dega a = k a$
and $\dega b = l b$. Now let $c\in \frac{[\omega]}{\nu}+
H^2_{\mbox{\rm\tiny dR}}(M)[[\nu]]$ be written as $c
=\frac{[\omega]}{\nu}+ \sum_{i=0}^\infty \nu^i c^i$. For the proof
of i.) we choose closed two-forms $\Omega_i$ such that
$\Omega_{2l+1} = 0$ (to achieve $[\Omega_{2l+1}]=c^{2l}=[0]$) and
$[\Omega_{2l+2}]= c^{2l+1}$ for all $l \in \mathbb N$ yielding $\Pa
\Omega = \Omega$. Moreover, we choose $s = \sum_{k=3}^\infty
s^{(k)}\in \W$ with $\sigma(s)=0$ and $\Pa s = s$. Under these
preconditions one easily proves that $\Pa r$ satisfies the
equations (\ref{Genr}) implying $\Pa r = r$ by uniqueness of the
solution of (\ref{Genr}). With such an element $r
\in \WL^1$ the Fedosov derivation $\mathcal D$ obviously commutes
with $\Pa$ implying that $\Pa \tau(f) = \tau (\Pa f)$ for all $f
\in C^\infty(M)[[\nu]]$ since $\Pa$ obviously commutes with
$\sigma$. Using this equation and the definition of $*$ together
with $\Pa ( (\Pa a) \circ (\Pa b))= (-1)^{k l} b \circ a$ and
observing that $\dega \tau(f)=0$ one gets the asserted property of
$*$ under the mapping $\Pa$. From Theorem \ref{Delclatheo} we get
$c(*)=\frac{[\omega]}{\nu} + \frac{1}{\nu} [\Omega] = c$. For ii.)
one proceeds quite analogously. The only difference lies in other
suitable choices of $\Omega$ and $s$, i.~e. we choose closed
two-forms $\Omega_i$ such that $\C \Omega_{2l+2} = \Omega_{2l+2}$,
$[\Omega_{2l+2}]= c^{2l+1}$ and $\C \Omega_{2l+1}= -
\Omega_{2l+1}$, $[\Omega_{2l+1}]= c^{2l}$ for all $l \in \mathbb
N$ implying $\C \Omega = \Omega$. Moreover, we choose $s\in \W$ such that
$\C s=s$. As in the proof of i.) one gets that $\C r = r$ yielding the
desired behaviour of the corresponding star product $*$ under the mapping
$\C$ as in the proof of part i.). The fact that $c(*) = c$ again follows
from Theorem \ref{Delclatheo} and the choice of $\Omega$. For the proof of
part iii.) one just has to bring into line the choices made for i.) with
the ones made for ii.) i.~e. choose $s$ with $\C s = s = \Pa s$ and closed
two-forms $\Omega_i$ with $\C \Omega_{2l+2}=\Omega_{2l+2}$ and
$\Omega_{2l+1}=0$ such that $[\Omega_{2l+2}] = c^{2l+1}$ and
$[\Omega_{2l+1}]=c^{2l}=[0]$ for all $l\in \mathbb N$. Then the argument as
in i.) and ii.) yields the stated result.
\end{PROOF}

\begin{REMARK}
The interest in such special star products from the viewpoint of
physics is based on the interpretation of the star product algebras
$(C^\infty(M)[[\nu]],\star)$ as the algebra of observables of the
quantised system corresponding to the classical system described by
the symplectic manifold $M$, and hence the existence of a
${}^*$-structure incorporated by complex conjugation (the
${}^*$-structure of the algebra of classical observables) is
strongly recommended. Moreover, the Weyl-Moyal product on
$T^*\mathbb R^n$ giving a correct description of the quantisation
of observables that are polynomials in the coordinates is of Weyl
type, motivating the general interest in such star products (cf.
\cite{BNW98,BNW99} for further details). In addition there is the
possibility of constructing ${}^*$-representations for star
products with $\C$ as ${}^*$-structure under the precondition of
having defined a formally positive functional on a suitable
twosided ideal in $C^\infty(M)[[\nu]]$ that is stable under $\C$ by
a formal analogue of the GNS construction (cf. \cite{BorWal98} for
details).
\end{REMARK}

To conclude this section we shall discuss the question of existence
of special equivalence transformations between equivalent star
products satisfying equations (\ref{Nupar}) and (\ref{Sternstr}).

\begin{DEFINITION}\label{SpeStarEquiDef}
Let $(C^\infty(M)[[\nu]],\star_1)$ and $(C^\infty(M)[[\nu]],
\star_2)$ denote equivalent star product algebras.
\begin{enumerate}
\item
In case $(C^\infty(M)[[\nu]],\star_1)$ and $(C^\infty(M)[[\nu]],
\star_2)$ have $\C$ incorporated as ${}^*$-structure they are
called equivalent as ${}^*$-algebras (resp. $\C$-equivalent) if
there is an equivalence transformation $\mathcal S$ between them
satisfying $\C \mathcal S \C = \mathcal S$.
\item
In case $(C^\infty(M)[[\nu]],\star_1)$ and
$(C^\infty(M)[[\nu]],\star_2)$ have the $\nu$-parity property they
are called $\Pa$-equivalent if there is an equivalence
transformation $\mathcal S$ between them satisfying $\Pa \mathcal S
\Pa = \mathcal S$.
\item
In case $(C^\infty(M)[[\nu]],\star_1)$ and $(C^\infty(M)[[\nu]],
\star_2)$ are of Weyl-type they are called Weyl-equiva\-lent if there
is an equivalence transformation $\mathcal S$ between them
satisfying $\C \mathcal S \C = \mathcal S$ and $\Pa
\mathcal S \Pa = \mathcal S$.
\end{enumerate}
\end{DEFINITION}

The following proposition states that for two equivalent star
products enjoying the additional algebraic properties as in the
preceding definition there are always equivalence transformations
being compatible with the mappings $\C$ and $\Pa$.

\begin{PROPOSITION}
Let $(C^\infty(M)[[\nu]],\star_1)$ and $(C^\infty(M)[[\nu]],
\star_2)$ denote equivalent star product algebras.
\begin{enumerate}
\item
In case $\C$ is a ${}^*$-structure for $\star_1$ and $\star_2$,
then $(C^\infty(M)[[\nu]],\star_1)$ and $(C^\infty(M)[[\nu]],
\star_2)$ are equivalent as ${}^*$-algebras.
\item
In case $\star_1$ and $\star_2$ have the $\nu$-parity property then
$(C^\infty(M)[[\nu]],\star_1)$ and $(C^\infty(M)[[\nu]],
\star_2)$ are $\Pa$-equivalent.
\item
In case $\star_1$ and $\star_2$ are of Weyl-type then $(C^\infty(M)
[[\nu]],\star_1)$ and $(C^\infty(M)[[\nu]], \star_2)$ are
Weyl-equivalent.
\end{enumerate}
\end{PROPOSITION}
\begin{PROOF}
For the proof of part i.) we consider some equivalence
transformation $\mathcal T$ between $\star_1$ and $\star_2$
satisfying $\mathcal T(f \star_1 g) = (\mathcal T
f)\star_2(\mathcal T g)$ for all $f,g
\in C^\infty(M)[[\nu]]$. Obviously $\C \mathcal T \C$ is also an
equivalence transformation between $\star_1$ and $\star_2$ and
hence there is an automorphism $\mathcal A$ of $\star_1$ such that
$\C \mathcal T \C = \mathcal T \mathcal A$. Conjugating this
equation with $\C$ and using $\C^2= \id$ we obtain $\mathcal T
=\C \mathcal T \C \C \mathcal A \C = \mathcal T
\mathcal A \C \mathcal A \C$ yielding $\mathcal A \C \mathcal A\C
=\id$. Since any automorphism of $\star_1$ has the shape
$\mathcal A = \exp(\nu \mathsf D)$ where $\mathsf D$ is a
derivation of $\star_1$ we get $\id =\exp(\nu \mathsf D) \exp(-\nu
\C \mathsf D \C)$ implying $\C \mathsf D \C = D$. For $t \in
\mathbb R$ we consider the automorphisms $\mathcal A^t :=
\exp( t \nu \mathsf D)$ of $\star_1$ satisfying $\C \mathcal A^t
\C = (\mathcal A^t)^{-1} = \mathcal A^{-t}$. Now $\mathcal S_t :=
\mathcal T \mathcal A^t$ obviously is an equivalence between
$\star_1$ and $\star_2$ for all $t \in \mathbb R$ and we have $\C
\mathcal S_t \C = \C \mathcal T \C \mathcal A^{-t} = \mathcal T
\mathcal A^{1-t} = \mathcal S_{1-t}$. Therefore $\mathcal S :=
\mathcal S_{1/2}$ satisfies $\mathcal S ( f \star_1 g) =
(\mathcal S f) \star_2 (\mathcal S g)$ and $\C \mathcal S \C =
\mathcal S$ proving part i.) of the proposition. For the proof of
part ii.) one proceeds completely analogously replacing $\C$ by
$\Pa$ in the above argumentation. For the proof of part iii.) we
consider some equivalence transformation $\mathcal T$ between $\star_1$
and $\star_2$ and use the results of part i.) and part ii.) to obtain two
further equivalence transformations $\mathcal S_1 = \mathcal T
\mathcal A_1^{1/2}$ and $\mathcal S_2 = \mathcal T \mathcal A_2^{1/2}$
satisfying $\C \mathcal S_1 \C = \mathcal S_1$ and $\Pa \mathcal S_2 \Pa =
\mathcal S_2$, where $\mathcal A_1$ and $\mathcal A_2$ are
automorphisms of $\star_1$ given by $\C \mathcal T \C = \mathcal T
\mathcal A_1$ and $\Pa \mathcal T \Pa = \mathcal T \mathcal A_2$.
In general $\mathcal S_1$ fails to satisfy $\Pa \mathcal S_1 \Pa =
\mathcal S_1$ as well as $\mathcal S_2$ fails to commute with $\C$,
but by an analogous procedure as for the proofs of the statements
i.) and ii.) $\mathcal S_1$ and $\mathcal S_2$ can be modified such
that the resulting equivalence transformations have the desired
properties. Since $\Pa$ commutes with $\C$ we have $\C \Pa \mathcal
T \Pa \C = \Pa \C
\mathcal T \C \Pa$ implying the crucial equation $\mathcal A_1 \C \mathcal A_2 \C =
\mathcal A_2 \Pa \mathcal A_1 \Pa$ by the definitions of $\mathcal A_1$ and
$\mathcal A_2$. Now we compute $\C \mathcal S_2\C = \mathcal S_2
\mathcal A_2^{-1/2} \mathcal A_1 \C \mathcal A_2^{1/2} \C=\mathcal S_2
\mathcal F_2$, where $\mathcal F_2:= \mathcal A_2^{-1/2} \mathcal A_1 \C
\mathcal A_2^{1/2} \C$ is an automorphism of $\star_1$ and hence $\mathcal
F_2=\exp(\nu \mathsf D_2)$ with a derivation $\mathsf D_2$ of $\star_1$.
As in ii.) one gets $\mathcal F_2 \C \mathcal F_2 \C = \id$ and $\mathcal R_2
:= \mathcal S_2 \mathcal F_2^{1/2}$ with $\mathcal F_2^{1/2}:=
\exp(\frac{\nu}{2}\mathsf D_2)$ is an equivalence transformation between
$\star_1$ and $\star_2$ satisfying $\C \mathcal R_2 \C = \mathcal R_2$. It
remains to show that $\mathcal R_2$ satisfies $\Pa \mathcal R_2 \Pa =
\mathcal R_2$. To this end we compute $\Pa \mathcal F_2 \Pa$ using
$\mathcal A_1 \C= \mathcal A_2 \Pa \mathcal A_1 \Pa \C \mathcal
A_2^{-1}$
\[
\Pa \mathcal F_2 \Pa= \Pa \mathcal A_2^{1/2}\Pa \mathcal A_1 \Pa \C \mathcal
A_2^{-1/2} \C \Pa = \mathcal A_2^{-1/2} \mathcal A_1 \C \mathcal A_2^{1/2}\C=
\mathcal F_2.
\]
Thus we find $\Pa \mathcal R_2 \Pa = \Pa \mathcal S_2 \Pa \Pa \mathcal F_2^{1/2}
\Pa= \mathcal S_2 \mathcal F_2^{1/2}=\mathcal R_2$ proving part iii.).
One can also modify $\mathcal S_1$ to obtain another equivalence transformation
$\mathcal R_1= \mathcal S_1\mathcal F_1^{1/2}$ having the desired properties
where $\mathcal F_1 :=\mathcal A_1^{-1/2} \mathcal A_2 \Pa
\mathcal A_1^{1/2} \Pa$ again is an automorphism of $\star_1$.
\end{PROOF}
\begin{REMARK}
The assertion about the existence of equivalence transformations
between equivalent star products with a ${}^*$-structure
incorporated by $\C$ that commute with $\C$ has an important
consequence for the GNS representations one can construct for these
star product algebras namely that such an equivalence
transformation induces a unitary map between the GNS Hilbert spaces
obtained by the GNS construction relating the corresponding GNS
representations (cf. \cite[Proposition 5.1]{BNW99}).
\end{REMARK}
\begin{appendix}
\section{The Deformed Cartan Formula}
\label{CarSec}
The aim of this section is to prove the deformed Cartan formula
that was very useful for our computations in Section
\ref{DelClaSec}. This formula and the proof of it which we shall
give already appeared in \cite[Lemma 4.6.]{Bor96}. A similar result
has also been derived in \cite[Proposition 4.3.]{Kra98} where the
vector field with respect to which the Lie derivative is computed
is assumed to be affine with respect to the symplectic connection
$\nabla$.
\begin{PROPOSITION}\label{LieDerProp}
For all vector fields $X \in \Gamma^\infty(TM)$ the Lie derivative
$\Lie{X}: \WL \to \WL$ can be expressed in the following manner:
\BEQ {LieDerallg}
    \Lie{X} = \mathcal D i_a(X) + i_a (X) \mathcal D + i_s(X) +
    dx^i \vee i_s( \nabla_{\partial_i} X) + \frac{1}{\nu}
    \ad(i_a(X)r).
\EEQ
In case $X=X_f$ is the Hamiltonian vector field of a function $f
\in C^\infty(M)$, i.~e. $i_{X_f}\omega = df$ this formula takes
the following form:
\BEQ {LieDerHam}
    \Lie{X_f} = \mathcal D i_a(X_f) + i_a (X_f) \mathcal D -
    \frac{1}{\nu} \ad \left( f + df \otimes 1 + \frac{1}{2}D df
    \otimes 1- i_a(X_f) r\right),
\EEQ
where $D = dx^i \vee \nabla_{\partial_i}$ denotes the operator of
symmetric covariant derivation.
\end{PROPOSITION}
\begin{PROOF}
The proof of formula (\ref{LieDerallg}) is obtained by collecting
the following formulas the proofs of which are all straightforward
computations just using the definitions of the involved mappings
and applying them to factorised sections $a = A \otimes \alpha
\in \WL$:
\BEQAR {LieDerallgbew}
    \label{1}\delta i_a(X) + i_a (X) \delta &=& i_s(X)\\
    \label{2}\frac{1}{\nu}\left(\ad (r) i_a(X) + i_a(X) \ad(r)
    \right) &=& \frac{1}{\nu} \ad ( i_a(X) r)\\
    \label{3}(\nabla i_a(X) + i_a(X) \nabla) (A \otimes \alpha)
    &=& \nabla_X  A \otimes \alpha + A \otimes \Lie{X} \alpha.
\EEQAR
For a symmetric one-form $A$ it is easy to see that $\nabla_X A=
\Lie{X} A - dx^i \vee i_s( \nabla_{\partial_i} X) A$. Together with
the observation that the operators on both sides of this equation
are derivations with respect to the $\vee$-product this and
(\ref{3}) imply
\BEQ {4}
    \nabla i_a(X) + i_a(X) \nabla = \Lie{X} - dx^i \vee i_s(
    \nabla_{\partial_i}X).
\EEQ
Combining (\ref{1}), (\ref{2}) and (\ref{4}) we get the first
statement of the proposition. For the second statement one just has
to observe that $\ad(f) = 0$ and that
\BEQAR {LiederHambew}
    i_s(X_f) = - \frac{1}{\nu} \ad ( df \otimes 1)&\quad& dx^i \vee
    i_s(\nabla_{\partial_i} X_f) = - \frac{1}{\nu} \ad\left(
    \frac{1}{2} D df \otimes 1\right),
\EEQAR
which is again a straightforward computation in local coordinates
using the explicit shape of the deformed product $\circ$. Using
these equations combined with (\ref{LieDerallg}) finishes the proof
of (\ref{LieDerHam}).
\end{PROOF}

\section{Computation of ${C_2^-}^\sharp$}
\label{CtwoSec}
This section just gives a sketch of the computations that are
necessary to determine the anti-symmetric part of the
bidifferential operator $C_2$ that occurs in the expression of the
Fedosov star product $f*g = fg + \nu C_1(f,g)+
\nu^2 C_2(f,g)+...$ of two functions $f,g
\in C^\infty(M)$.
\begin{PROPOSITION}\label{CtwoProp}
The anti-symmetric part $C_2^-$ of the bidifferential operator
$C_2$ is given by
\BEQ {CtwoGleich}
    C_2^-(f,g)= \frac{1}{2}\left( C_2(f,g) - C_2(g,f)\right) = -
    \frac{1}{2} \left( \Omega_1 + ds_1^{(3)}\right)(X_f, X_g)=
    \rho_2(X_f,X_g),
\EEQ
where $f,g \in C^\infty(M)$ and $X_f$ resp. $X_g$ denote the
corresponding Hamiltonian vector fields with respect to $\omega$
and $s_1^{(3)} \in \Gamma^\infty(T^*M)$ denotes the one-form
occurring in the first order of $\nu$ in $s^{(3)}=s_3^{(3)}+ \nu
s_1^{(3)}$, where $s_3^{(3)}\in \Gamma^\infty(\bigvee^3 T^*M)$,
that comes up from the normalisation condition $\delta^{-1} r = s$
(cf. equation (\ref{Genr})). Thus we have
\BEQ {Cwoclass}
    c(*)^0 = -2{C_2^-}^\sharp = -2 [\rho_2] = [\Omega_1].
\EEQ
\end{PROPOSITION}
\begin{PROOF}
Using the shape of the fibrewise product $\circ$ we obtain
$f*g-g*f= \nu \sigma (\Lambda^{rs}i_s(\partial_r) \tau(f)
i_s(\partial_s) \tau(g)) + O(\nu^3)$. To compute the terms of order
less or equal to two in $\nu$ we thus only have to know $\tau(f)$
and $\tau(g)$ except for terms of symmetric degree and $\nu$-degree
greater than one. Hence it is enough to look at
$\tau(f)^{(0)},\ldots,\tau(f)^{(3)}$ since for $\tau(f)^{(k)}$ with
$k\geq 4$ either the symmetric degree or the $\nu$-degree of the
occurring terms are greater than one. Looking at the recursion
formula (\ref{GenTaylorRecurs}) we thus see that the only terms of
$r$ that are needed are given by $r^{(2)}= \delta s^{(3)}$ and
$r^{(3)}= \delta s^{(4)} +
\delta^{-1}(R + \nu \Omega_1 + \nabla r^{(2)} - \frac{1}{\nu}
r^{(2)} \circ r^{(2)} )$ what is obtained from (\ref{GenrRecus}) by
writing down the terms of total degree $2$ resp. $3$. Writing
$\approx$ for equations holding modulo terms of symmetric degree
resp. $\nu$-degree greater than one, one gets by lengthy but
obvious computation that
\begin{eqnarray*}
\tau(f)^{(0)}&=& f\\
\tau(f)^{(1)}&=& df \otimes 1\\
\tau(f)^{(2)}&=&\frac{1}{2}Ddf\otimes 1 - i_s(X_f)s_3^{(3)}
\approx 0\\
\tau(f)^{(3)}&\approx& -\delta^{-1}(i_s(X_f) r^{(3)})\approx
-\nu \left( i_s(X_f)s_2^{(4)} + \frac{1}{2}i_{X_f}\left(\Omega_1 +
ds_1^{(3)}\right) \right)\otimes 1,
\end{eqnarray*}
where we have written $s^{(4)}= s_4^{(4)} + \nu s_2^{(4)}$ with
$\degs s_k^{(4)} = ks_k^{(4)}$. Inserting these results into
$f*g-g*f$ as given above the terms involving $s_2^{(4)} \in
\Gamma^\infty(\bigvee^2 T^*M)$ cancel because of their symmetry and
one gets
\[
f*g-g*f = \nu \{f,g\} - \nu^2 (\Omega_1 + ds_1^{(3)})(X_f,X_g) +
O(\nu^3)
\]
proving the proposition.
\end{PROOF}

One should observe, that this is the only instance of our proof of
Theorem \ref{Delclatheo} where the modified normalisation condition
on $r$ enters our considerations, whereas the other terms of $c(*)$
could be computed without making use of it.
\end{appendix}

\section*{Acknowledgements}
\label{AcknoSec}
I would like to thank Martin Bordemann for many useful discussions
and suggestions, in particular for pointing out reference
\cite{Bor96} containing the deformed Cartan formula. Moreover, I
would like to thank Stefan Waldmann for many valuable discussions
via e-mail during his stay at the Math Department of UC Berkeley.
Finally I should like to thank Alexander V. Karabegov for the
inspiring discussion during his visit at the seminar on fundamental
interactions at the University of Freiburg.

\end{document}